\documentclass{amsart}
\usepackage{amsrefs,mathrsfs}
\usepackage[xnumber]{math}
\usepackage{lineno}\linenumbers
\def\pair<#1>{{\langle\!\langle}#1{\rangle\!\rangle}}

\newtheorem{mainthm}{Theorem}
\newtheorem{maincor}[mainthm]{Corollary}

\begin{document}
\title{(Self-)similar groups and the Farrell-Jones conjectures}
\author{Laurent Bartholdi}
\date{26 July 2011}
\address{Mathematisches Institut\\Georg-August Universit\"at zu
  G\"ottingen\\ Bunsenstra\ss e 3--5\\ D-37073 G\"ottingen\\ Germany}
\email{laurent.bartholdi@gmail.com}
\thanks{The work is supported by the Courant Research Centre ``Higher
  Order Structures'' of the University of G\"ottingen}
\begin{abstract}
  We show that contracting self-similar groups satisfy the Farrell-Jones
  conjectures as soon as their universal contracting cover is
  non-positively curved. This applies in particular to bounded
  self-similar groups.

  We define, along the way, a general notion of contraction for groups
  acting on a rooted tree in a not necessarily self-similar manner.
\end{abstract}
\maketitle

\section{Introduction}
Few properties are known to hold for all groups; in the recent years,
counterexamples have been found to numerous ``plausible conjectures'',
usually formulated as questions: is there an infinite, finitely
generated group all of whose elements have finite order? is there an
amenable group that cannot be produced using extensions and filtered
colimits of virtually abelian groups? is there a group whose word
growth is strictly between polynomial and exponential?

The ``Farrell-Jones conjectures'', predicting how the algebraic
K-/L-theory of the group ring $RG$ may be expressed in terms of the
algebraic K-/L-theory of $R$ and the group theory of $G$, is one of
the prominent remaining
conjectures~\cite{bartels-luck-reich:fjappl}. If it is satisfied by
the group $G$, numerous group-theoretical consequences for $G$ follow,
in particular $RG$ has no non-trivial idempotent if $G$ is
torsion-free and $R$ is a domain of characteristic $0$. The
Farrell-Jones conjectures are inherited under many group-theoretical
operations (finite direct and free products, filtered colimits), but
possibly not under wreath products; we say the Farrell-Jones
conjectures \emph{hold with wreathing} if they hold for all wreath
products $G\wr P$ with a finite permutation group $P$.

In search of a possible counterexample to the Farrell-Jones
conjectures, it might have been speculated that the ``self-similar
groups'' studied by Al\"eshin, Grigorchuk, Gupta and Sidki since the
1970s would play an important role; indeed, these groups have served
to answer or illuminate all the questions in the first paragraph.

Self-similar groups are groups acting in a recursive manner on a
regular rooted tree $T_d$. If the recursion of every element involves
only a linearly growing subtree of $T_d$, the group is said to
be \emph{bounded}.

We show in this note that considerable care will be required to
construct a counterexample within the class of self-similar groups. We
prove (see below for precise definitions):
\begin{mainthm}
  Let $G$ be a bounded self-similar group. Then $G$ satisfies the
  Farrell-Jones conjectures.
\end{mainthm}

\begin{mainthm}
  Let $G$ be a contracting similar group. Then $G$ satisfies the
  Farrell-Jones conjectures if its universal contracting cover satisfies
  the Farrell-Jones conjectures with wreathing.
\end{mainthm}

\begin{maincor}
  The Al\"eshin-Grigorchuk, Gupta-Sidki, GGS, and generalized
  Grigorchuk groups all satisfy the Farrell-Jones conjectures.
\end{maincor}

\subsection{Acknowledgments}
Wolfgang L\"uck encouraged me to write this short note, with the
intent of narrowing the domains of group theory in which a
counterexample is to be searched.

Thomas Schick and Wolfgang L\"uck generously provided valuable
feedback on a preliminary version, and clarified for me the status of
the Farrell-Jones conjectures with respect to wreath products.

\section{The Farrell-Jones conjectures}
We review very briefly the statement of the Farrell-Jones conjectures;
we include them for definiteness, but will never work directly with
their definition.

A model for the \emph{virtually cyclic classifying space} $E^{vc}(G)$
is a topological $G$-space $X$ whose isotropy groups are all virtually
cyclic, and such that for any topological $G$-space $Y$ with virtually
cyclic isotropy groups there exists up to $G$-homotopy a unique
$G$-map $Y\to X$.

The Farrell-Jones conjectures assert that the natural map
\[H_n^G(E^{vc}(G),\mathbf S)\to H_n^G(\{.\},\mathbf S),\] induced by
$E^{vc}(G)\to\{.\}$, is an Farrell-Jones for all $n$. Here $\mathbf
S$ is either the K-theory spectrum $\mathbf K_{\mathcal A}$ or the
L-theory spectrum $\mathbf L_{\mathcal A}^{\langle-\infty\rangle}$
over the orbit category associated with an additive $G$-category
$\mathcal A$.

For our purposes, it suffices to note that the class of groups for
which the conjectures are known to hold contains virtually abelian
groups, hyperbolic groups~\cite{bartels-luck:hypcat0} for $n\le1$,
CAT(0) groups~\cites{bartels-luck:hypcat0,wegner:cat0}, cocompact
lattices in virtually connected Lie groups, threefold
groups~\cite{bartels-luck:lattices} and arithmetic groups over
algebraic number fields (unpublished). It is closed under taking
subgroups, colimits~\cite{bartels-luck:crossed}*{Corollary~0.8}, and
finite direct and free products. (This is the advantage of using the
more general version with co\"efficients in an additive category ---
the inheritance properties come almost for free).

Note that, in general, it is not known whether the conjectures are
inherited under finite extensions. Since every finite extension is a
subgroup of the wreath product with a finite
group~\cite{kaloujnine-krasner:extensions}, the question reduces to
whether the conjecture is inherited by finite wreath products. This is
known in some specific cases, in particular for cocompact lattices in
virtually connected Lie groups, threefold groups, arithmetic groups
over algebraic number fields, and CAT(0) groups, as we now explain.

CAT(0) spaces are metric spaces in which triangles are at least as
thin as in euclidean space; see the classical
reference~\cite{bridson-h:msnpc}.  CAT(0) groups, also called
\emph{non-positively curved groups}, are groups acting properly,
isometrically and cocompactly on a CAT(0) space of finite topological
dimension. That class contains virtually abelian groups, and is closed
under direct, free and finite wreath products.

\begin{lem}\label{lem:cat0}
  If $G$ is CAT(0), then so is $G\wr P$ for any finite permutation
  group $P$.
\end{lem}
\begin{proof}
  Let $G$ act properly discontinuously on the CAT(0) space $X$, and
  let $P$ be a permutation group on $n$ points. Then $G\wr
  P=G^n\rtimes P$ acts properly discontinuously on $X^n$, with $G^n$
  acting co\"ordinatewise and $P$ by permutation of the co\"ordinates.
\end{proof}

\section{(Self-)similar groups}
We summarize the notion of \emph{self-similar group}, presenting it in
a slightly more general and algebraic manner than is usual;
see~\cite{nekrashevych:ssg} or~\cite{bartholdi-g-s:bg} for classical
references. By $G\wr d$ we denote the permutational wreath product
$G^d\rtimes\sym{d}$.

A \emph{self-similar group} is a group $G$ endowed with a homomorphism
$\phi:G\to G\wr d$, called its \emph{self-similarity structure}. The
integer $d$ is the \emph{degree} of the self-similarity
structure. Usually, the self-similarity is implicit, and one simply
denotes by $G$ the self-similar group.

The map $\phi$ can be applied diagonally to all entries in $G^d$,
yielding a map $G^d\to(G\wr d)^d$, and therefore a map $G\wr d\to(G\wr
d)\wr d\subseteq G\wr(d^2)$; more generally, we get maps $G\wr d^n\to
G\wr d^{n+1}$ which we all denote by $\phi$. We may compose these
maps, and write $\phi^n$ for the iterate $\phi^n:G\to G\wr d^n$.

By projecting to the permutation part, we then have homomorphisms
$G\to\sym{d^n}$ and, assembling them together, a permutational action
of $G$ on $T_d:=\bigsqcup_{n\ge0}\{1,\dots,d\}^n$; one may identify
$T_d$ with the vertex set of a rooted $d$-regular tree, by connecting
$v_1\dots v_n$ to $v_1\dots v_nv_{n+1}$ for all $v_i\in\{1,\dots,d\}$,
in such a way that $G$ acts by graph isometries. This action need not
be faithful; if it is, then $G$ is called a \emph{faithful
  self-similar group}.

A self-similar group is \emph{contracting} if there exists a finite
subset $N\subset G$ such that, for all $g\in G$ and all $n$ large
enough, $\phi^n(g)\in N^{d^n}\times\sym{d^n}$. The
smallest such $N$ is called the \emph{nucleus} of $G$.

Let $\tilde F$ denote the free group on $N$. By definition, the
nucleus satisfies the condition $\phi(N)\subset
N^d\times\sym{d}$. The restriction of $\phi$ to $N$ can
therefore uniquely be extended to a homomorphism $\tilde\phi:\tilde
F\to\tilde F\wr d$. Set
\[R=\{w\in N\cup N^2\cup N^3\subset\tilde F\mid w=_G1\}.\]
Similarly, we have $\tilde\phi(R)\subset R^d\times1$. Set $F=\tilde
F/R$. The homomorphism $\tilde\phi$ then induces a homomorphism, again
written $\phi:F\to F\wr d$.

Note that $F$ is a finitely presented group, and that the natural map
$N\subset\tilde F\to N\subset G$ defines a homomorphism $F\to G$. We
will see in Lemma~\ref{lem:contr} that $F$ is contracting, with
nucleus $N$. However, the self-similarity structure of $F$ need not be
faithful, even if that of $G$ was faithful. We call $F$ the
\emph{universal contracting cover} of $G$. Note also that in general
the homomorphism $F\to G$ need not be onto, or equivalently $N$ need
not generate $G$. This is, however, the case in all examples we
present here.

Here are some extreme examples; more classical ones appear
in~\S\ref{s:examples}. The full group $W$ of isometries of $T_d$ is
self-similar, but not contracting; actually not even countable. Its
subgroup $\{g\in W\mid \psi^n(g)\in\{1\}^{d^n}\times\sym{d^n}\text{
  for some }n\}$ is faithful, self-similar, and contracting with
nucleus $\{1\}$. Any group $G$, with $\phi:G\to G^d$ the diagonal
embedding, defines a non-faithful self-similar structure on $G$, which
is contracting precisely when $G$ is finite. Consider finally $A$ a
finite group, and $G$ the group of finitely-supported functions $\Z\to
A$. Take $d=2$, and set $\psi(f)=\pair<f_0,f_1>$ with $f_0(n)=f(2n)$
and $f_1(n)=f(2n-1)$. This defines a self-similarity structure on $G$,
which is not faithful, and contracting with nucleus
$N=\{\text{functions supported on }\{0,1\}\}$. Our main result does
not give any interesting information on such actions.

\subsection{Similar groups}
We now generalize the definitions above to more general groups. A
group $G$ is \emph{similar} if there exists a sequence
$G=G_0,G_1,\dots$ of groups, a sequence of integers $d_1,d_2,\dots$,
and a sequence of homomorphisms $\phi_n:G_n\to G_{n+1}\wr
d_{n+1}$. The similarity structure is faithful if the corresponding
permutational action on
$\bigsqcup_{n\ge0}\{1,\dots,d_1\}\times\cdots\times\{1,\dots,d_n\}$ is
faithful. Again abusing notation, the compositions of $\phi_n$'s are
written $\phi_n^m:G_n\to G_{n+m}\wr d_{n+1}d_{n+2}\cdots d_{n+m}$.

Let $N_0,N_1,\dots$ be a sequence of finite sets, with $N_n\subset
G_n$ for all $n$. We say that $G$ \emph{contracts to $(N_n)_{n\ge0}$}
if for every $g\in G_n$ and every $m$ large enough, $\phi_n^m(g)\in
N_{n+m}^{d_{n+1}d_{n+2}\cdots d_{n+m}}\times\sym{d_{n+1}d_{n+2}\cdots
  d_{n+m}}$.

In that case, it is possible, up to enlarging the $N_n$'s, to assume
$\phi_n(N_n)\subset
N_{n+1}^{d_{n+1}}\times\sym{d_{n+1}}$, and we always
make that additional assumption. We call the sequence $N_0,N_1,\dots$
a \emph{nucleus} of $G$.

Note however that the sequence $N_0,N_1,\dots$ is not unique --- for
example, it is always possible to replace finitely many of the initial
terms by $1$. We say $G$ is \emph{generated by its nucleus} if $N_n$
generates $G_n$ for all $n$.

Extending the previous definition, let $F_n$ be the finitely presented
group
\[F_n:=\langle N_n\mid\text{ words of length $\le3$ that are $\equiv1$
  in }G_n\rangle.
\]
We then have induced homomorphisms $F_n\to F_{n+1}\wr d_{n+1}$,
defining a similarity structure for the group $F:=F_0$.

\begin{lem}\label{lem:contr}
  The similar group $F$ contracts to $(N_n)_{n\ge0}$.
\end{lem}
\noindent We again call $F$ the \emph{universal contracting cover} of
$G$; note that it depends on the choice of $(N_n)_{n\ge0}$.
\begin{proof}
  Consider $n\in\N$. For every $g\in N_n^{\le2}\subset G_n$, there
  exists $m\in\N$ such that $\phi_n^m(g)\in
  N_{n+m}^{d_{n+1}d_{n+2}\cdots
    d_{n+m}}\times\sym{d_{n+1}d_{n+2}\cdots d_{n+m}}$,
  by the contraction condition. Since there are finitely many $g$'s
  under consideration, there exists $m_n\in\N$ such that
  \[\phi_n^{m_n}(N_n^2)\in N_{n+m_n}^{d_{n+1}d_{n+2}\cdots
    d_{n+m_n}}\times\sym{d_{n+1}d_{n+2}\cdots d_{n+m_n}}.
  \]
  On the other hand, consider $\tilde w\in\tilde F_n$ a word of length
  $\ell\le2$ in the alphabet $N_n$, and denote by $\overline w$ and
  $w$ respectively its image in $G_n$ and in $F_n$. The entries in
  $\tilde\phi_n^{m_n}(\tilde w)$ have length precisely $\ell$, by
  construction. They are termwise equal, in $G_{n+m_n}$, to the
  entries of $\phi_n^{m_n}(\overline w)$. Since $F_{n+m_n}$ contains
  all relations of length $\le3$, these entries are also termwise
  equal in $F_{n+m_n}$. It follows that, for every $w\in F_n$ of
  length $\le2$, all entries of $\phi_n^{m_n}(w)$ all belong to
  $N_{n+m_n}$.

  Consider now $g\in F_n$, of length $\ell\le 2^k$ in the alphabet
  $N_n$. Set inductively $n_0=n$ and $n_{i+1}=n_i+m_{n_i}$. By the
  previous paragraph, the entries of $\phi_n^{m_n}(g)$ have length
  $\le2^{k-1}$ over $N_{n+m_n}=N_{n_1}$, and more generally the
  entries of $\phi_n^{n_k-n}$ have length $\le2^0$ in $N_{n_k}$, that
  is, they belong to $N_{n_k}$.
\end{proof}

We call a similar group \emph{contracting} if it has been endowed with
a sequence $(N_n)_{n\ge0}$ to which it contracts. Note that this fixes
the choice of a contracting finitely presented cover. Similar
contracting groups naturally include self-similar groups, by
considering constant sequences $G$, $\phi$, $N$ and $F$.

Note that we explicitly allow the sequences $G,\phi$ to be constant
while the $N_n$'s increase. Quite generally, if each $G_n$ is
countable, then there exists a sequence of finite sets to which it
contracts; namely, enumerate $G_n=\{g_{n,1},g_{n,2},\dots\}$, and let
$N_n$ be the set of co\"ordinates of $\phi_m^{n-m}(g_{m,i})$ for all
$i,m\le n$. Understandably, our main result applies formally to such
constructions, but does not yield any useful information.

\subsection{Main result}
\begin{prop}\label{prop:main}
  Let $G$ be a faithful contracting similar group, generated by its
  nucleus. If all terms $F_n$ of the universal contracting cover of
  $G$ satisfy the Farrell-Jones conjectures with wreathing, then $G$
  satisfies the Farrell-Jones conjectures.
\end{prop}
\begin{proof}
  In the self-similar case, set $K_0=1\triangleleft F$, and
  $K_{n+1}=\phi^{-1}(K_n^d)$ for all $n\ge0$; and finally
  $K_\infty=\bigcup_{n\ge0}K_n$.  More generally, in the similar case,
  set $K_n=\ker(\phi^n)\triangleleft F$ and
  $K_\infty=\bigcup_{n\ge0}K_n$.

  There is an natural homomorphism $\pi:F/K_\infty\to G$, which we
  prove to be an Farrell-Jones. Let $g\in F$ be in the kernel of $\pi$;
  then, because $F$ is contracting, there is $n\in\N$ such that
  $\phi^n(g)$ belongs to $N_n^{d_1\cdots
    d_n}\times\sym{d_1\cdots d_n}$; furthermore, the
  permutation is trivial because $\phi^n\pi(g)=\phi^n(1)=1$, and the
  entries in $N_n$ are trivial because $F_n$ contains relations of
  length $1$ in $N_n$. Therefore $g\in K_n$ so $g\in K_\infty$, as was
  to be shown.

  We then have $G=\lim F/K_n$, and because the Farrell-Jones conjectures
  are stable under colimits it suffices to see that $F/K_n$ satisfies
  the Farrell-Jones conjectures. By the first Farrell-Jones theorem,
  $F/K_n$ is a subgroup of $F_n\wr d_1\cdots d_n$, so it suffices to
  show that $F_n\wr d_1\cdots d_n$ satisfies the Farrell-Jones
  conjectures. Since $F_n$ satisfies the Farrell-Jones conjectures with
  wreathing, we are done.
\end{proof}

As stated in the introduction, Proposition~\ref{prop:main} applies in
particular to contracting similar groups whose universal contracting
cover are CAT(0) groups, lattices in virtually connected Lie groups,
or arithmetic groups over algebraic function fields.

\section{Examples}\label{s:examples}
We now give some examples of contracting, similar groups, recall some
of their basic properties, and show that they satisfy the
Farrell-Jones conjectures.

We follow a slightly unorthodox path to define (self-)similar groups:
we first give their contracting covers, and then simply say that the
group itself is the faithful quotient of the cover. This, of course,
defines uniquely the self-similar group $G$ in question: it is the
quotient of its universal contracting cover $F$ by the normal subgroup
$K_\infty\triangleleft F$.

We denote by $\pair<g_1,\dots,g_d>\sigma$ an element of the wreath
product $G\wr d$, with $\sigma$ written as a product of disjoint
cycles.

\subsection{The Al\"eshin and Grigorchuk groups}
The Al\"eshin-Grigorchuk group is obtained as follows. Set
\[F=\langle a,b,c,d\mid a^2,b^2,c^2,d^2,bcd\rangle=C_2*(C_2\times C_2),\]
and define $\phi:F\to F\wr2$ by
\[\phi(a)=\pair<1,1>(1,2),\quad\phi(b)=\pair<a,c>,\quad\phi(c)=\pair<a,d>,\quad\phi(d)=\pair<1,b>.\]
Let $G$ be the faithful self-similar quotient of $F$.

This group (up to finite index) was first considered
in~\cite{aleshin:burnside}, providing a ``tangible'' example of
infinite, finitely generated, torsion group (the first examples of
groups with these properties are due to
Golod~\cite{MR28:5082}). Grigorchuk proved in~\cite{grigorchuk:growth}
that its word growth is strictly between polynomial and exponential,
and in~\cite{grigorchuk:amenEG} that it is amenable, but not
elementary amenable. It is contracting, with nucleus $\{1,a,b,c,d\}$.

Since $F$ is CAT(0), as a free product of finite groups, $G$ satisfies
the Farrell-Jones conjectures by Proposition~\ref{prop:main}.

More elaborate examples have also been constructed by
Grigorchuk~\cite{grigorchuk:gdegree}. Fix an infinite sequence
$\omega=\omega_0\omega_1\cdots$ of epimorphisms $(C_2\times
C_2)\cong\langle b,c,d\rangle\to\langle a\rangle\cong C_2$, and assume
that $\omega$ contains infinitely many of each of the three possible
epimorphisms.  Define homomorphisms $\phi_n:F\to F\wr2$ for all
$n\ge0$ by
\[\phi(a)=\pair<1,1>(1,2),\quad\phi(x)=\pair<\omega_n(x),x>\text{ for }x\in\{b,c,d\}.\]
Let $G_\omega$ be the faithful similar quotient of $F$ using this
similarity structure.

Again, $G_\omega$ is contracting with nucleus $N_n=\{1,a,b,c,d\}$ for
all $n\in\N$, so all such groups satisfy the Farrell-Jones
conjectures. There are uncountably many such groups, and they all are
torsion $2$-groups of intermediate word growth.

\subsection{The Gupta-Sidki groups}
The Gupta-Sidki groups are obtained as follows. Choose a prime $p\ge3$,
set
\[F=\langle a,t\mid a^p,t^p\rangle=C_p*C_p,\]
and define $\phi:F\to F\wr p$ by
\[\phi(a)=\pair<1,\dots,1>(1,\dots,p),\quad\phi(t)=\pair<a,a^{-1},1,\dots,1,t>.\]
Let $G$ be the faithful self-similar quotient of $F$.

These groups are shown in~\cite{gupta-s:burnside} to be infinite,
finitely-generated torsion $p$-groups.

Since $F$ is CAT(0), as a free product of finite groups, $G$ satisfies
the Farrell-Jones conjectures by Proposition~\ref{prop:main}.

\subsection{Bounded groups}
Assume that $G$ is a self-similar group, and that, for every $g\in G$,
there exists a bound $B\in\N$ such that, for all $n\in\N$, there are
at most $B$ non-trivial entries in $\phi^n(g)$. Note that it suffices
to check this property for the generators of $G$; and that it holds
for the generators of the Grigorchuk group with $B=2$, and those of
the Gupta-Sidki groups for $B=3$.

It is then known (see~\cite{bondarenko-n:pcf}) that $G$ is
contracting. More precisely, $G$ is isomorphic to a subgroup of a
self-similar group of very special type
(see~\cite{bartholdi-k-n-v:ba}). Fix an integer $d\ge2$, set
\[F=\sym{d}*(\sym{d}\wr\sym{d-1}),\]
and define $\phi:F\to F\wr d$ by
\[\phi(\sigma)=\pair<1,\dots,1>\sigma,\quad\phi(g:=\pair<f_1,\dots,f_{d-1}>\tau)=\pair<f_1,\dots,f_{d-1},g>\tau.\]

Since $F$ is CAT(0), as a free product of finite groups, $G$ satisfies
the Farrell-Jones conjectures by Proposition~\ref{prop:main}.

Note that the faithful quotient of $F$ is amenable; this is
how~\cite{bartholdi-k-n-v:ba} show that all bounded self-similar
groups are amenable.

\subsection{Dynamics}
Let $f$ be a branched covering of a topological space $\mathcal M$;
this means there is an open dense subset $\mathcal M_0\subset\mathcal
M$ and a covering $f:\mathcal M_0\to\mathcal M$. We assume $f$ has
finite degree $d$. Let $P_f$ denote the \emph{post-critical locus} of
$f$:
\[P_f=\bigcup_{n\ge1}f^n(\mathcal M\setminus \mathcal M_0).
\]
Assume finally that $\mathcal M\setminus P_f$ is
path-connected. Choose a basepoint $*$, and for each $x\in f^{-1}(*)$
choose an arc $\ell_x$ from $*$ to $x$ in $\mathcal M\setminus
P_f$. Number also $f^{-1}(*)$ as $\{x_1,\dots,x_d\}$.

These data define a self-similar group as follows. It is again defined
via a cover, $F:=\pi_1(\mathcal M\setminus P_f,*)$. Consider
$\gamma\in F$. For each $x_i\in f^{-1}(*)$, let $\gamma_i$ denote the
unique $f$-lift of $\gamma$ that starts at $x_i$, and let it end at
$x_{\pi(i)}\in f^{-1}(*)$. Define then $\phi:F\to F\wr d$ by
\[\phi(\gamma)=\pair<\ell_{\pi(1)}^{-1}\gamma_1\ell_1,\dots,\ell_{\pi(d)}^{-1}\gamma_d\ell_d>\pi.\]

If $\mathcal M$ is in fact a locally simply connected metric space,
and $f$ is uniformly expanding (meaning there exists $\lambda>1$ such
that $d(fx,fy)>\lambda d(x,y)$ whenever $d(x,y)$ is sufficiently
small), then $F$ is contracting.

This applies in particular to $\mathcal M$ a complex manifold and $f$
a holomorphic map (which is then expanding for the Kobayashi metric).

The special case $\mathcal M=\mathbb C$ and $f$ a degree-$2$
polynomial has been extensively studied
in~\cite{bartholdi-n:mandelbrot1}. The cover $F$ is a free group, so
this provides more examples of groups satisfying the Farrell-Jones
conjectures. One important such example, associated with the map
$f(z)=z^2-1$, has been studied in~\cite{MR2003h:60011}
and~\cite{bartholdi-v:amenability}; it is amenable, orderable, of
exponential growth, and residually poly-$\mathbb Z$.

Other examples, on higher-dimensional manifolds, have been considered
by Koch et al.~\cites{koch:phd,buff-e-k-p:pullback}. There, the
universal contracting cover is the sphere braid group.

\section{Conclusion}
We have shown that if a counter-example to the Farrell-Jones conjectures
exists in the class of (self-)similar groups, it will not be an easy
matter to establish that fact.

For one thing, with very few exceptions, non-contracting self-similar
groups are intractable (it required considerable effort to prove that
the elementary example of~\cite{aleshin:free} is a free group!)

For another, calculations in a contracting self-similar groups are
usually reduced to calculations in a finitely presented group, in
which one may manipulate words. It would be surprising that the
Farrell-Jones conjectures fail for self-similar group, yet be unsettled
for its cover.

Since the Farrell-Jones conjectures is not settled for the sphere braid
group, we have, at the present, no argument to check the Farrell-Jones
conjectures on the faithful self-similar quotient of the braid groups
that arise in this manner.

Let $G$ be a self-similar group, and let $e\in\N$ be given. Assume
that, for every $g\in G$, there exists a bound $B\in\N$ such that, for
all $n\in\N$, there are at most $Bn^e$ non-trivial entries in
$\phi^n(g)$. Then $G$ is said to be of \emph{polynomial activity
  growth} of degree $e$; see~\cite{sidki:polynonfree}, who proves that
such groups do not contain free subgroups.

It is then known~\cite{amir-a-v:linamen} that $G$ embeds, possibly for
larger $d$, in a specific group $P(d,e)$ of polynomial activity
growth, defined by its cover as follows. Set $\Sigma_{-1}=\sym{d}$ and
$\Sigma_i=\Sigma_{i-1}\wr\sym{d-1}$ for $i=0,\dots,e$; set
\[F=\Sigma_{-1}*\cdots*\Sigma_e,\]
and define $\phi:F\to F\wr d$ by
\begin{align*}
  \phi(\sigma) &= \pair<1,\dots,1>\sigma,\\
  \phi(g) &= \pair<f_1,\dots,f_{d-1},g>\tau\text{ for }g=\pair<f_1,\dots,f_{d-1}>\tau\in\Sigma_i,i\ge0.
\end{align*}
These are non-contracting self-similar groups if $e\ge1$; for $e\le1$,
the faithful quotient is
amenable~\cites{bartholdi-v:amenability,amir-a-v:linamen}, while
amenability of the faithful quotient is open for larger $e$.

The arguments in~\cite{sidki:polynonfree} show that the nucleus $N$ of
$P(d,e)$, while infinite, admits a partial well ordering, such that
every $g\in N$ has the form $g\in\Sigma_{-1}$ or
$\phi(g)=\pair<g_1,\dots,g_{d-1},g>$ with $g_i<g$ for all
$i\in\{1,\dots,d-1\}$. Presumably this means that arguments similar to
those given here show that $P(d,e)$, and therefore all its subgroups,
satisfy the Farrell-Jones conjectures.

It has been conjectured by Nekrashevych that all contracting
self-similar groups are amenable; although no conclusive link has been
established between amenability and the Farrell-Jones conjectures.

At the other extreme of contracting self-similar groups lie
\emph{bireversible groups}. These are self-similar groups $(G,\phi)$
such that the map $G\times\{1,\dots,d\}\to G\times\{1,\dots,d\}$,
given by $(g,i)\mapsto(g_i,\sigma(i))$ if
$\phi(g)=\pair<g_1,\dots,g_d>\sigma$, is a bijection. They are related
to the infinite simple groups constructed
in~\cites{burger-m:localglobal,burger-m:fpsimple}. They would seem
like a natural class in which to look at counterexamples, though all
examples studied up to now are lattices in virtually connected Lie
groups.

\end{document}